\documentclass[12pt,twoside]{article}
\usepackage[hmargin=3cm,vmargin=3cm]{geometry}
\usepackage[T1]{fontenc} 
\usepackage{mathptmx}    
\usepackage[latin1]{inputenc}
\usepackage{amsmath,amsthm,amssymb,amsfonts}
\usepackage{enumerate}
\usepackage{graphicx,psfrag}
\usepackage{url}
\usepackage{pst-node,pstricks-add,pst-text,pst-3d,pst-tree}

\newcommand{\reff}[1]{(\ref{#1})}

\newcommand{\mysection}[1]{\paragraph{\uppercase{\bf #1}.}}

\theoremstyle{plain}
\newtheorem{theorem}{Theorem}

\newtheorem*{theorem*}{Theorem}
\newtheorem*{lemma*}{Lemma}
\newtheorem*{fact*}{Fact}
\newtheorem*{corollary*}{Corollary}

\theoremstyle{definition}

\begin{document}

\begin{center}
\LARGE \bf A Dilution Test for the Convergence of \\ Subseries of a Monotone Series
\end{center}

\begin{center}
\large \bf Lasse Leskel\"a and Mikko Stenlund
\end{center}

\vspace{5mm}

\mysection{1. Introduction}

Cauchy's condensation test allows to determine the convergence of a monotone series by looking at a
weighted subseries that only involves terms of the original series indexed by the powers of two. It
is natural to ask whether the converse is also true: \emph{Is it possible to determine the
convergence of an arbitrary subseries
\[
 \sum_{k \ge 1} a_{s(k)} = a_{s(1)} + a_{s(2)} + a_{s(3)} + \cdots
\]
of a monotone series $\sum_{n \ge 1} a_n$ by looking at a suitably weighted version of the original
series?} In this note we show that the answer is affirmative and introduce a new convergence test
particularly designed for this purpose.

\mysection{2. Cauchy's and Schl\"omilch's condensation tests}

Consider a series which is monotone in the sense that its terms satisfy $a_1 \ge a_2 \ge \cdots \ge
0$. Cauchy's condensation test (e.g.\ \cite[Theorem 2.3]{Bonar_Khoury_2006}) states that a monotone
series $\sum_{n \ge 1} a_n$ converges if and only if
\[
 \sum_{k \ge 0} 2^k a_{2^k} =
 a_1 + 2 a_2 + 4 a_4 + 16 a_{16} + \cdots
\]
converges, thereby allowing to determine the convergence of a monotone series by only looking at
its terms indexed by the powers of two. Schl\"omilch's extension (e.g.\ \cite[Theorem
2.4]{Bonar_Khoury_2006}) allows to replace the powers of two by a more general subsequence $s(1) <
s(2) < s(3) < \cdots$ of the positive integers, assuming that the forward differences
\begin{equation}
 \label{eq:ForwardDifferences}
 \Delta s(k) = s(k+1) - s(k)
\end{equation}
do not grow too fast.

\begin{theorem}[Schl\"omilch]
\label{the:Condensation}
For any monotone series $\sum_{n \ge 1} a_n$ and subsequence of the integers such that for
some $c>0$,
\begin{equation}
 \label{eq:ForwardBound}
 \frac{\Delta s(k+1)}{\Delta s(k)} \le c \quad \text{for all $k \ge 1$},
\end{equation}
the series $\sum_{n \ge 1} a_n$ converges if and only if $\sum_{k \ge 1}  a_{s(k)} \, \Delta s(k)$
converges.
\end{theorem}

Cauchy's condensation test can be recovered as a special case of Theorem~\ref{the:Condensation} by
substituting $s(k) = 2^{k-1}$. We will next present a short proof of Schl\"omilch's result to
highlight its structural similarity to our new test given in Section~3.

\begin{proof}[Proof of Theorem~\ref{the:Condensation}]
Because the terms of $\sum_{n \ge 1} a_n$ are nonincreasing, we see that
\begin{equation}
 \label{eq:BasicInequality}
 a_{s(k+1)} \Delta s(k)
 \ \le \ a_{s(k)} + \cdots + a_{s(k+1)-1}
 \ \le \ a_{s(k)} \Delta s(k)
\end{equation}
for all $k \ge 1$. Inequalities~\reff{eq:BasicInequality} combined with
assumption~\reff{eq:ForwardBound} imply that
\[
 c^{-1} a_{s(k+1)} \Delta s(k+1)
 \ \le \ a_{s(k)} + \cdots + a_{s(k+1)-1}
 \ \le \ a_{s(k)} \Delta s(k).
\]
By summing the above display over $k$ we now find that
\[
 c^{-1} \sum_{k \ge 2} a_{s(k)} \Delta s(k)
 \ \le \ \sum_{n \ge s(1)} a_n
 \ \le \ \sum_{k \ge 1} a_{s(k)} \Delta s(k),
\]
so that all three series above either converge or diverge together.
\end{proof}

\mysection{3. Dilution test}

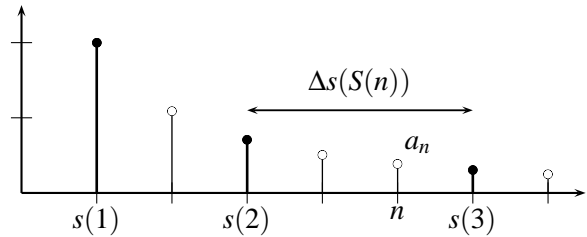
\begin{figure}[h]
\begin{center}
\psset{unit=1cm}
\begin{pspicture}(-0.2,-.2)(7.5,2.5)
  \psaxes[labels=none]{->}(7.5,2.5)

  \psline[linewidth=0.5pt](1,0)(1,2)
  \psline[linewidth=0.5pt](2,0)(2,1.0887)
  \psline[linewidth=0.5pt](3,0)(3,0.7071)
  \psline[linewidth=0.5pt](4,0)(4,0.5060)
  \psline[linewidth=0.5pt](5,0)(5,0.3849)
  \psline[linewidth=0.5pt](6,0)(6,0.3054)
  \psline[linewidth=0.5pt](7,0)(7,0.2500)
  \psdots*[dotscale=1, dotstyle=o](1,2)(2,1.0887)(3,0.7071)(4,0.5060)(5,0.3849)(6,0.3054)(7,0.2500)

  \psline[linewidth=1pt](1,0)(1,2)
  \psline[linewidth=1pt](3,0)(3,0.7071)
  \psline[linewidth=1pt](6,0)(6,0.3054)
  \psdots*[dotscale=1](1,2)(3,0.7071)(6,0.3054)

  \psline[linewidth=.7pt]{<->}(3,1.1)(6,1.1)
  \uput[90](4.5,1.1){\small $\Delta s(S(n))$}

  \uput[-90](5,0){\small $n$}
  \uput[50](5,0.3849){\small $a_n$}

  \uput[-90](1,0){\small $s(1)$}
  \uput[-90](3,0){\small $s(2)$}
  \uput[-90](6,0){\small $s(3)$}

\end{pspicture}
\end{center}
\caption{\label{fig:Gap} The weight $\Delta s(S(n))$.}
\end{figure}

Cauchy's and Schl\"omilch's condensation tests are designed for determining the convergence of a
monotone series by looking at a weighted subseries of the original series. We will now reverse this
line of thought and prove a converse to these results, which allows to determine the convergence of
a subseries of a monotone series using a weighted version of the original series. Suitable weights
can be defined in terms of the forward differences~\reff{eq:ForwardDifferences} and the counting
function
\[
 S(n) = \# \{k: s(k) \le n\}
\]
of a sequence $s(1) < s(2) < \cdots$. It is interesting to note that the growth condition of the
forward differences in Theorem~\ref{the:Condensation} is not needed below. The weights $\Delta
s(S(n))$ in Theorem~\ref{the:Dilution} measure the distance between points of the subsequence
nearest to $n$; see Figure~\ref{fig:Gap}.

\begin{theorem}
\label{the:Dilution}
For any monotone series $\sum_{n \ge 1} a_n$ and any infinite subsequence of the integers, the
subseries $\sum_{k \ge 1} a_{s(k)}$ converges if and only if
\begin{equation}
 \label{eq:Dilution}
 \sum_{n \ge s(1)} \frac{a_{n}}{\Delta s(S(n))} < \infty.
\end{equation}
\end{theorem}

\begin{proof}
The monotonicity of the series implies the validity of~\reff{eq:BasicInequality}. After dividing
the terms in~\reff{eq:BasicInequality} by $\Delta s(k)$, we find that
\begin{equation}
 \label{eq:Inequality2}
 a_{s(k+1)}
 \ \le \ \frac{a_{s(k)}}{\Delta s(k)} + \cdots + \frac{a_{s(k+1)-1}}{\Delta s(k)}
 \ \le \ a_{s(k)}.
\end{equation}
Because $S(n)$ indexes the last member of the sequence $s(1) < s(2) < \cdots$ not exceeding $n$, it
follows that $S(n) = k$ for all $n$ such that $s(k) \le n \le s(k+1)-1$. This is
why~\reff{eq:Inequality2} may be rephrased as
\[
 a_{s(k+1)}
 \ \le \! \sum_{n=s(k)}^{s(k+1)-1} \frac{a_n}{\Delta s(S(n))}
 \ \le \ a_{s(k)}.
\]
After summing the above display over $k$, we find that
\[
 \sum_{k \ge 2} a_{s(k)}
 \ \le \sum_{n \ge s(1)} \frac{a_n}{\Delta s(S(n))}
 \ \le \ \sum_{k \ge 1} a_{s(k)},
\]
which shows that all series above either converge or diverge together.
\end{proof}

\mysection{4. Thinning out a divergent series into a convergent one}

Given a divergent monotone series $\sum_{n \ge 1} a_n$, one may ask whether it can be made
convergent by deleting some of its terms. If $\lim_{n \to \infty} a_n > 0$, this is obviously not
possible, while if $\lim_{n \to \infty} a_n = 0$ this can always be done by selecting terms of the
series along a sparse enough subsequence. Indeed, in this case the series may even be thinned out
to sum to an arbitrary positive real number (Banerjee and Lahiri~\cite{Banerjee_Lahiri_1964}).

For a divergent monotone series $\sum_{n \ge 1} a_n$ such that $a_n \to 0$, a more specific
question is to quantify a sufficient degree or sparsity required for the thinning subsequence. The
following corollary of the dilution test presents a sufficient condition.

\begin{theorem}
\label{the:Sparsity}
Consider a monotone divergent series $\sum_{n \ge 1} a_n$ such that $\sum_{n \ge 1} a_n^p$
converges for some $p>1$. A sufficient condition for the convergence of the subseries $\sum_{k \ge
1} a_{s(k)}$ is that
\begin{equation}
 \label{eq:Sparsity}
 \sum_{k \ge 1} s(k)^{-1/p} < \infty.
\end{equation}
\end{theorem}
\begin{proof}
Observe that the series~\reff{eq:Dilution} in Theorem~\ref{the:Dilution} can be written as $\sum_{k
\ge 1} A_k$, where
\[
 A_k = \frac{1}{\Delta s(k)} \sum_{n=s(k)}^{s(k+1)-1} a_n
\]
is the average of the terms $a_{s(k)}, \dots, a_{s(k+1)-1}$. Because all these terms are less than
or equal to the terms $a_1, \dots, a_{s(k)}$, we see that $A_k$ is bounded from above by the
average
\[
 B_k = \frac{1}{s(k)} \sum_{n=1}^{s(k)} a_n.
\]
Jensen's inequality implies that
\[
 B_k^p
 \le \frac{1}{s(k)} \sum_{n=1}^{s(k)} a_n^p
 \le \frac{1}{s(k)} \sum_{n=1}^{\infty} a_n^p,
\]
so that
\[
 \sum_{n \ge s(1)} \frac{a_{n}}{\Delta s(S(n))}
 = \sum_{k \ge 1} A_k
 \le  \sum_{k \ge 1} B_k
 \le \left( \sum_{n=1}^{\infty} a_n^p \right)^{1/p} \sum_{k \ge 1} s(k)^{-1/p}.
\]
The now claim follows as a consequence of Theorem~\ref{the:Dilution}.
\end{proof}

\mysection{5. Sparse subseries of the harmonic series}

The harmonic series $\sum_{n \ge 1} \frac{1}{n}$ is probably the best-known example of a divergent
series; see \cite[Section 3]{Bonar_Khoury_2006} for a lively discussion.
Kempner has shown~\cite{Kempner_1914} that, rather surprisingly, we obtain a convergent series
by deleting from the harmonic series all terms whose decimal representation contains the digit `9'.
Kempner's curious series has afterwards attracted lots of interest, with several articles
generalizing and sharpening the basic result; see for example Schmelzer and
Baillie~\cite{Schmelzer_Baillie_2008} and references therein.

The following corollary of Theorem~\ref{the:Sparsity} shows that the harmonic series converges over
any polynomially sparse subsequence. A subsequence of the integers is called \emph{polynomially
sparse} if its density among the first $n$ positive integers decreases fast enough as $n$ grows,
according to
\begin{equation}
 \label{eq:PolynomiallySparse}
 S(n)/n \le c n^{-\alpha}
\end{equation}
for some $c>0$ and $\alpha \in (0,1)$. A simple counting argument (e.g.
Behforooz~\cite{Behforooz_1995}) may be used to verify that Kempner's no-`9'
fulfills~\reff{eq:PolynomiallySparse} with $c=10$ and $\alpha = 1-\frac{\log 9}{\log 10}$.


\begin{theorem}
\label{the:Harmonic}
The harmonic series converges over any polynomially sparse subsequence.
\end{theorem}
\begin{proof}
Fix an integer $k \ge 1$, and let $c$ and $\alpha$ be such that \reff{eq:PolynomiallySparse} holds
for all $n$. The definition of the counting function implies that $s(k) \ge n+1$ for all integers
$n$ such that $S(n) < k$, and in particular for all integers $n$ such that $n < (k/c)^\beta$, where
$\beta = 1/(1-\alpha)$. By letting $n$ be the largest integer strictly less than $(k/c)^\beta$, we
see that $s(k) \ge (k/c)^\beta$. Therefore, condition~\reff{eq:Sparsity} of
Theorem~\ref{the:Sparsity} is valid for any $p \in (1,\beta)$. The claim now follows by
Theorem~\ref{the:Sparsity}, because the series $\sum_{n \ge 1} n^{-p}$ converges for all $p>1$
(this well-known fact is usually proved by using Cauchy's condensation test).
\end{proof}


Theorem~\ref{the:Harmonic} may also be proved as a consequence of a stronger result specialized to
the harmonic series (Powell and \v{S}al{\'a}t \cite{Powell_Salat_1991}): The harmonic series over a
subsequence converges if and only if the counting function of the subsequence satisfies $\sum_{n
\ge 1} S(n)/n^2 < \infty$.

\mysection{6. Concluding remark}

Many nonmonotone series $\sum_{n\geq 1}a_n$ encountered in applications admit a monotone majorant series $\sum_{n\geq 1}b_n$. 
In this case, the dilution test can be applied to subseries of the majorant series; if $\sum_{k\geq 1} b_{s(k)} $ converges, then so does the corresponding subseries $\sum_{k\geq 1}a_{s(k)}$ of the original nonmonotone series.




\vspace{2ex}

{\small \noindent {\bf ACKNOWLEDGEMENTS.} Both authors have been supported by the Academy of
Finland. Mikko Stenlund is also affiliated with the University of Helsinki.}

{\small
\bibliographystyle{abbrv}
\bibliography{lslReferences}
}

\vspace{1ex}

{\small

\noindent
{\it Aalto University, PO Box 11100, 00076 Aalto, Finland \\
lasse.leskela@iki.fi}

\vspace{.5ex}

\noindent
{\it Courant Institute, New York, NY 10012-1185, USA \\
mikko@cims.nyu.edu}

}

\end{document}